%
%
%
\input amstex
\magnification=1200
\documentstyle{amsppt}
\topmatter
\title
Criterion for the resolvent set of nonsymmetric tridiagonal operators
\endtitle
\author
A.I. Aptekarev, V. Kaliaguine and W. Van Assche
\endauthor
\affil
Keldysh Institute of Applied Mathematics \\
Nizhni\u\i\ Novgorod State University \\
Katholieke Universiteit Leuven
\endaffil
\address
Keldysh Institute of Applied Mathematics,
Russian Academy of Science,
Miuss\-ka\-ya Sq.\ 4,
125047 Moscow, Russia
\endaddress
\email
aptekaa\@applmat.msk.su
\endemail
\address
Department of Mathematics,
Nizhni\u\i\ Novgorod State University,
Gagarina 23A,
Nizhni\u\i\ Novgorod, Russia
\endaddress
\email
mechmat\@nnucnit.nnov.su
\endemail
\address
Department of Mathematics,
Katholieke Universiteit Leuven,
Celestijnenlaan 200\,B,
B-3001 Heverlee, Belgium
\endaddress
\email
fgaee03\@cc1.kuleuven.ac.be
\endemail
\thanks
The first author is partly supported by the Russian Basic Research Foundation
(grant 93-01278); the third author is
a Senior Research Associate of the Belgian National Fund
for Scientific Research
\endthanks
\keywords
tridiagonal operators, resolvents, Pad\'e approximation
\endkeywords
\subjclass
47A10, 41A21
\endsubjclass
\leftheadtext{A.I. Aptekarev, V. Kaliaguine and W. Van Assche}
\rightheadtext{nonsymmetric tridiagonal operators}

\abstract
We study nonsymmetric tridiagonal operators acting in the Hilbert space
$\ell^2$ and describe the spectrum and the resolvent set of such
operators in terms of a continued fraction related to the resolvent.
In this way we establish a connection between Pad\'e approximants
and spectral properties of nonsymmetric tridiagonal
operators.
\endabstract

\endtopmatter

\document
\head 1. Introduction \endhead
In this paper we consider operators in the Hilbert
space $\ell^2$, with the following representation in an orthonormal basis
of this space:
$$ A=\pmatrix
\beta_0 & \gamma_0 & 0 & 0 & \ldots \\
\alpha_1 & \beta_1 & \gamma_1 & 0 & \ldots \\
0 & \alpha_2 &\beta_2 & \gamma_2 & \ldots \\
\vdots & \vdots & \ddots & \ddots & \ddots \\
\endpmatrix   .         \tag 1 $$
where $\alpha_k$, $\beta_k$, $\gamma_k$ are  complex numbers,
$\alpha_k \neq 0, \gamma_k \neq 0, k \in {\Bbb N}$.
The operator $A$ is defined for any
finite vector $x=x_0g_0+x_1g_1+ \ldots +x_ng_n$
in the orthonormal basis $\{g_n\}_0^{\infty}$, and its domain of
definition $D(A)$ is dense in $\ell^2$.

The symmetric case, where $\alpha_k=\overline{\gamma_{k-1}}$ and $\beta_k$
real, can under
an appropriate chosen basis be reduced to an infinite Jacobi matrix
enabling a deeper examination of this spectral theory. The main tools of
investigation
in this case are the classical moment problem (including the theory of
selfadjoint extensions of unbounded symmetric operators), the theory of general
orthogonal polynomials and the spectral theorem for selfadjoint operators
(\cite{1}, \cite{2}, \cite{9}, \cite{11}). This
connection between the
spectral theory and analysis is fruitful for various points of view. For
example, the scattering problem for a Jacobi matrix can be treated on the basis
of strong (or Szeg\H{o} type) asymptotic results for orthogonal polynomials
(\cite{4}, \cite{6}, \cite{10}). On the other hand,
the perturbation theory gives new results for orthogonal polynomials
(\cite{12}).

In the general, nonsymmetric, case we can not use the spectral theorem
and orthogonal polynomials, but in return we will obtain results of use
for Pad\'e approximants and the theory of continued fractions.
We denote as usual by $\sigma(A)$ the spectrum of the
operator $A$, $\Omega(A) = {\Bbb C} \setminus \sigma(A)$ is the resolvent set
and the resolvent is
$R(\lambda)=(\lambda I-A)^{-1}$
for $\lambda \in \Omega$. By $\{g_n\}_0^{\infty}$ we denote the orthonormal
basis in $\ell^2$. The function
$$ \phi(\lambda)=(R(\lambda)g_0,g_0)    \tag 2 $$
is analytic on the resolvent set $\Omega(A)$. If the operator $A$ is bounded
(this is the case when $\alpha_k,\beta_k, \gamma_k$ are bounded),
then the function
$\phi(\lambda)$ is analytic for $|\lambda| > \| A \| $.
For convenience we have changed the basis of representation
of the operator $A$. The
new basis $\{e_n\}_0^{\infty}$ is defined by $e_n=g_n/d_n$, where $d_n=$
$\gamma_0\gamma_1 \ldots \gamma_{n-1}, d_0=1$. In this basis the operator
\thetag{1} has the following form:
$$ A=\pmatrix
\beta_0 & 1 & 0 & 0 & \ldots \\
\alpha_1\gamma_0 & \beta_1 & 1 & 0 & \ldots \\
0 & \alpha_2\gamma_1 &\beta_2 & 1 & \ldots \\
\vdots & \vdots & \ddots & \ddots & \ddots \\
\endpmatrix . \tag 3 $$
If we set $a_n=\alpha_n\gamma_{n-1}$, $b_n=\beta_n$, then we have in the
basis $\{e_n\}$:
$$ A=\pmatrix
b_0 & 1 & 0 & 0 & \ldots \\
a_1 & b_1 & 1 & 0 & \ldots \\
0 & a_2 & b_2 & 1 & \ldots \\
\vdots & \vdots & \ddots & \ddots & \ddots \\
\endpmatrix .   \tag 4 $$
The Chebyshev algorithm
(Algorithm 7.2.1 in \cite{8, page 248})
 applied to $\phi(\lambda)$ at infinity gives us
the following continued fraction (J-fraction):
$$ \frac{1|}{|\lambda-b_0}-\frac{a_1|}{|\lambda-b_1}-\frac{a_2|}{|\lambda-b_2}-
\cdots\  . \tag 5 $$
The numerators $P_n$ and denominators $Q_n$ of the $n$th convergents for this
fraction satisfy the three-term recurrence relation:
$$ a_ny_{n-1}+b_ny_n+y_{n+1}=\lambda y_n, \qquad n=0,1,2,\ldots , \tag 6 $$
with the initial conditions
$$ \left\{ \matrix
Q_{-1}=0,\quad & Q_0=1, \\
P_{-1}=1,\quad & P_0=0.
\endmatrix \right. $$

The main problem considered in this paper is to describe the spectrum and
the resolvent
set of operators of type \thetag{1} in terms of the continued fraction
\thetag{5} and the monic polynomials $P_n, Q_n$ and in this way to
establish a connection between Pad\'e approximants and spectral properties of
nonsymmetric tridiagonal operators. In our main result, {\bf theorem 1}, we
state the criterion for the resolvent set $\Omega(A)$ in terms of
the growth of the polynomials
$Q_n(z)$.
The proof of the necessary condition is essentially
based on the results of Kershaw and Demko, Moss and Smith about
decay rates
of inverses of banded matrices (\cite{5}). This theorem gives us analytic
properties of the polynomials $Q_n(\lambda)$ and the related Pad\'e approximants
$\pi_n:=P_n/Q_n$ on the set $\Omega(A)$. We prove, among other
things, the following
\proclaim{Theorem 2}
If $A$ is bounded and
$$ 0 <  |\alpha_k|,|\gamma_k| \leq C_1 < \infty $$
for some $C_1$, then for any $\lambda \in \Omega(A)$ there is a
subsequence of Pad\'e approximants $\pi_n$ which converges to
$\phi(\lambda)$ with a geometric rate.
\endproclaim

In section 2 we formulate and prove theorem 1 and related results.
Section 3 will be devoted to an application of these results to the
behavior of Pad\'e approximants.

\head 2. Analysis of the spectrum and the resolvent set of the operator
\endhead
Suppose
$A$ has a representation \thetag{4} in an orthogonal basis $\{e_n\}$,
where $\| e_n \| = h_n$.
We assume that the operator $A$ is bounded.
In this case $A$ is defined on the whole space $\ell^2$ and $\lambda \in
\Omega(A)$ if and only if the following conditions are satisfied:
\roster
\item $Ker(\lambda I-A)=\{0\}$,
\item $\forall n,\exists z_n: (\lambda I-A)z_n=e_n$,
\item $\| (\lambda I-A)^{-1}x \| \leq C \| x \|,$
for all finite vectors $x=x_0e_0+x_1e_1+ \cdots +x_ne_n$ and
some constant $C$.
\endroster

\proclaim{Proposition 1}
$\lambda_0$ is an eigenvalue of $A$ if and only if
$$ \sum_{n=0}^{\infty}|Q_n(\lambda_0)|^2h_n^2 < +\infty .   \tag 7 $$
\endproclaim
\demo{Proof}
 We can formally write $Ax=\lambda_0x$ and obtain for
$x=x_0e_0+x_1e_1+ \cdots$ the following system
$$ \left\{     \matrix
  & & b_0x_0& + & x_1 & = & \lambda_0x_0 \\
  a_1x_0 & + & b_1x_1 & + & x_2 & = & \lambda_0x_1 \\
  a_2x_1 & + & b_2x_2 & + & x_3 & = & \lambda_0x_2 \\
   &  & \vdots & & & \vdots & \vdots
\endmatrix \right.     $$
If $x_0=0$ then $x_n=0, \forall n \in {\Bbb N}$. If $x_0 \neq 0$, then
$x_n=Q_n(\lambda_0)x_0$ and in this case $x \in \ell^2$ iff
$\sum_0^{\infty}|Q_n(\lambda_0)|^2h_n^2 < \infty$. \qed
\enddemo

\proclaim{Proposition 2}
For given $\lambda$, the equation
$$ (\lambda I-A)z=e_n   $$
has a solution in $\ell^2$ for all $n \geq 0$ if and only if there exists
a complex number $\gamma$ such that:
$$ \sum_{n=0}^{\infty}|Q_n(\lambda) \gamma -P_n(\lambda)|^2h_n^2 < +\infty .
        \tag 8 $$
\endproclaim
\demo{Proof}
 First we find the formal solution $z^{(0)}$ of the equation
$$ (\lambda I-A)z=e_0.  $$
Let $z^{(0)}=(z_{0,0}, z_{1,0}, z_{2,0}, \ldots)$ in the basis
$\{e_n\}$ then
$$ \left\{   \matrix
  1 & + & b_0z_{0,0} & + & z_{1,0} & = & \lambda z_{0,0} \\
  a_1z_{0,0} & + & b_1z_{1,0} & + & z_{2,0} & = & \lambda z_{1,0} \\
  a_2z_{1,0} & + & b_2z_{2,0} & + & z_{3,0} & = & \lambda z_{2,0} \\
   &  & \vdots & &  & \vdots &\vdots
   \endmatrix   \right. \tag 9 $$
So $z_{n,0}$ satisfies the same recurrence relations as $Q_n(\lambda)$
and $P_n(\lambda)$ for $n \geq 1$, only the first relation is different.
This implies that $z_{n,0}=Q_n(\lambda) \gamma -P_n(\lambda)$ for $n>1$
and for any $\gamma$ the first relation is also satisfied:
$b_0 \gamma + \gamma (\lambda-b_0)-1+1=\lambda \gamma.$ Thus the sequence
$Q_n(\lambda) \gamma -P_n(\lambda)$ is a solution of \thetag{9} for
any $\gamma$. We need a solution in $\ell^2$. So $e_0$ is in the image of
$(\lambda I-A)$ iff
$\sum_0^{\infty}|Q_n(\lambda) \gamma -P_n(\lambda)|^2h_n^2 < +\infty$
for some $\gamma$. Note that if $\lambda$ is not an eigenvalue of $A$ then
the convergence of this series is possible only for one $\gamma$.
For the solution of the equation $(\lambda I-A)z=e_1$ we note that
$(\lambda I-A)e_0= \lambda e_0-Ae_0= (\lambda -b_0)e_0-a_1e_1$
and then $-a_1e_1= (\lambda I-A)(e_0-(\lambda-b_0)z^{(0)})$. So
$z^{(1)}=-(e_0-(\lambda-b_0)z^{(0)})/a_1$ is a solution of
$(\lambda I-A)z=e_1$. This is an element of $\ell^2$ if $z^{(0)} \in \ell^2$.
In the same way we find $z^{(2)}=-(e_1-(\lambda-b_1)z^{(1)}+z^{(0)})/a_2$
for the solution of $(\lambda I-A)z=e_2$ and so on. \qed
\enddemo

\proclaim{Theorem 1}
Suppose $A$ is bounded, $\lambda$ is not an eigenvalue of $A$ and
\thetag{8}
is satisfied for some $\gamma$. Then
$\lambda \in \Omega(A)$ if and only if there are constants $C>0$ and
 $0<q<1$
such that
$$ |\frac{Q_n(\lambda)r_m(\lambda)}{a_1a_2\ldots a_n} \frac{h_m}{h_n}|,
 \leq  Cq^{m-n}, \qquad n \leq m, \tag 10 $$
$$ |\frac{Q_m(\lambda)r_n(\lambda)}{a_1a_2\ldots a_n} \frac{h_m}{h_n}|
 \leq  Cq^{n-m} , \qquad n \geq m, \tag 11 $$
where $r_m(\lambda):=Q_m(\lambda)\gamma -P_m(\lambda)$.
\endproclaim
\demo{Proof}
From the propositions 1 and 2 it follows that the inverse
operator $B=(\lambda I-A)^{-1}$ can be defined on the basis vectors by
$Be_n=z^{(n)}$. The matrix of $B$ in the basis $\{e_n\}$ is of the form:
$$ \pmatrix
r_0 & \frac{r_1}{a_1} & \frac{r_2}{a_1a_2} & \frac{r_3}{a_1a_2a_3} &
\ldots\\
r_1 & Q_1\frac{r_1}{a_1} & Q_1\frac{r_2}{a_1a_2} & Q_1\frac{r_3}{a_1a_2a_3}
& \ldots \\
r_2 & r_2\frac{Q_1}{a_1} & Q_2\frac{r_2}{a_1a_2} & Q_2\frac{r_3}{a_1a_2a_3}
& \ldots \\
r_3 & r_3\frac{Q_1}{a_1} & r_3\frac{Q_2}{a_1a_2} & Q_3\frac{r_3}{a_1a_2a_3}
& \ldots \\
\vdots & \vdots & \vdots & \vdots & \ddots
\endpmatrix .  \tag 12 $$
Indeed, for the first column there is no problem, since it contains
$z^{(0)}$. The
$n$th column contains $z^{(n)}$. For $z^{(1)}=-(e_0-(\lambda-b_0)z^{(0)})/a_1$
we have
$$ \aligned
z_{0,1} &= -\frac{1-(\lambda-b_0)r_0}{a_1}= \frac{r_1}{a_1}  \\
z_{1,1} &= -\frac{-(\lambda-b_0)z_{1,0}}{a_1} =Q_1 \frac{r_1}{a_1} \\
z_{2,1} &= -\frac{-(\lambda-b_0)z_{2,0}}{a_1} =r_2 \frac{Q_1}{a_1} \\
\vdots\ &\ \vdots \qquad \vdots
\endaligned \tag 13 $$
This is the second column of $B$. Suppose we have calculated $z^{(k)}$ as
indicated and we have to find $z^{(k+1)}$. First note that
$$
z^{(k+1)}=\frac{e_k-(\lambda-b_k)z^{(k)}+z^{(k-1)}}{a_{k+1}}  .
$$
Then
$$ \align
z_{0,k+1} &= -\frac{1}{a_{k+1}}
               (-(\lambda-b_k)\frac{r_k}{a_1a_2\cdots a_k} +
               \frac{r_{k-1}}{a_1a_2\cdots a_{k-1}}) =
               \frac{r_{k+1}}{a_1a_2\cdots a_{k+1}} \\
z_{1,k+1} &= -\frac{1}{a_{k+1}}
               (-(\lambda-b_k)\frac{Q_1r_k}{a_1a_2\cdots a_k} +
               \frac{Q_1r_{k-1}}{a_1a_2\cdots a_{k-1}}) =
               Q_1\frac{r_{k+1}}{a_1a_2\cdots a_{k+1}}  \\
\vdots &\vdots  \qquad \vdots  \\
z_{k,k+1} &= -\frac{1}{a_{k+1}}
               (1-(\lambda-b_k)\frac{Q_kr_k}{a_1a_2\cdots a_k} +
               \frac{Q_{k-1}r_k}{a_1a_2\cdots a_{k-1}}) =
               Q_k\frac{r_{k+1}}{a_1a_2\cdots a_{k+1}} .
\endalign $$
We used here the relations $Q_{k-1}r_k-Q_kr_{k-1}=-a_1a_2\cdots a_{k-1}$
deduced from the recurrence relations for $Q_k$ and $r_k$.
Next,
$$
z_{k+1,k+1} = -\frac{1}{a_{k+1}}
               (-(\lambda-b_k)\frac{Q_kr_{k+1}}{a_1a_2\cdots a_k} +
               \frac{Q_{k-1}r_{k+1}}{a_1a_2\cdots a_{k-1}}) =
               r_{k+1}\frac{Q_{k+1}}{a_1a_2\cdots a_{k+1}}
$$
and so on for $z_{l,k+1}$, $l>k+1$. Thus the form of the matrix $B$ is correct.
Now we decompose $B$ into two parts: $B=B_1+B_2$, where $B_1$ and $B_2$
are lower and upper triangular matrices which the same diagonal
elements $(B_1)_{i,i} = (B_2)_{i,i} = \frac12 B_{i,i}$.
For the estimation of norms of matrices
we use the formula
$$ \| K \| = \sup_{\| x \| \leq 1}
                        \sup_{\| y \| \leq 1}
                        |(y, Kx)| . $$
We have for two finite vectors $x=x_0e_0+x_1e_1+\cdots+x_ne_n$ and
$y=y_0e_0+y_1e_1+\cdots+y_ne_n$:
$$ \align
(y,B_1x) = &\frac{1}{2}\overline{b_{0,0}}y_0\overline{x_0}h_0^2 + \\
           &(\overline{b_{1,0}}y_1\overline{x_0} +\frac{1}{2}
    \overline{b_{1,1}}y_1\overline{x_1})h_1^2 +\cdots  \\
   & + (\overline{b_{n,0}}y_n\overline{x_0}+\overline{b_{n,1}}y_n\overline{x_1}
+ \cdots +\frac{1}{2}\overline{b_{n,1}}y_n\overline{x_n})h_n^2 .
\endalign $$
Here, we denote by $b_{i,j}$ the elements of the matrix $B$.
Then
$$ \align
|(y,B_1x)| &\leq  \frac{1}{2}\sum_{k=0}^n|b_{k,k}||x_k||y_k|h_k^2 +
 \sum_{j=1}^n \sum_{k=0}^{n-j}|b_{k+j,k}||x_k||y_{k+j}|h_{k+j}^2  \\
 &=  \frac{1}{2}\sum_{k=0}^n|b_{k,k}||x_k||y_k|h_k^2 +
 \sum_{j=1}^n
\sum_{k=0}^{n-j}|b_{k+j,k}\frac{h_{k+j}}{h_k}||x_kh_k||y_{k+j}h_{k+j}|  \\
 &\leq  C( \frac12 +q+q^2+ \cdots +q^n)\| x \|\ \| y \|\\
 &\leq  C\frac{1-q^{n+1}}{1-q}\| x \|\ \| y \| .
\endalign $$
A similar estimation is applicable for the matrix $B_2$. The  first part of
the theorem is thus proved.

To prove that for $\lambda \in \Omega(A)$  the estimation \thetag{10}
holds, one can use the known theorem on the decay rates of the inverse of banded
matrices (see \cite{5}) and the matrix representation of
$B=(\lambda I-A)^{-1}$
obtained in the first part of our proof. \qed
\enddemo

\head 3. Application to Pad\'e approximants
\endhead
In this section we consider some applications of theorem 1. First of all we
note  that theorem 1 together with propositions 1 and 2 gives a
characterization of the resolvent set of $A$ in terms of the polynomials
$P_n$ and $Q_n$. On the other hand this theorem also gives a result on the
convergence of the continued fractions \thetag{5}. We recall that the
fraction $P_n/Q_n$ is a diagonal type Pad\'e approximant for $\phi(\lambda)$
at infinity (see, e.g., Theorem 7.15 (B) in \cite{8, page 250}).
\proclaim{Corollary 2}
If $\lambda \in \Omega(A)$, then the remainder of Pad\'e approximation
in linear form tends to zero, i.e.,
$$  \lim_{n \to \infty} [Q_n(\lambda)\phi(\lambda) - P_n(\lambda)] h_n = 0 . $$
\endproclaim
\demo{Proof}
We have $z^{(0)} = R(\lambda) e_0$ and $\phi(\lambda) = ( R(\lambda)e_0,e_0 )
= z_{0,0}$, with $h_0=1$. On the other hand $z_{0,0} = \gamma Q_0 - P_0 =
\gamma$, hence $\phi(\lambda) = \gamma$ and the corollary follows from
proposition 2. \qed
\enddemo

\noindent{\bf Remark 1:} In the case of a symmetric operator $A$ with a
representation \thetag{1} in some orthogonal basis, the polynomials
$$    q_n = Q_n/d_n, \quad d_n = \gamma_0\gamma_1\cdots \gamma_{n-1} $$
are orthonormal and if $\lambda \in \Omega(A)$ then
$q_n(\lambda)\phi(\lambda)-p_n(\lambda) \to 0$ as $n \to \infty$. In this
case we have $h_n = 1/|d_n|$.

\proclaim{Corollary 3}
If there exists a positive constant $C_1$  such that
$|\alpha_k|, |\gamma_k| \leq C_1$
then for all $\lambda \in \Omega(A)$
$$   \limsup_{n \to \infty} |Q_n(\lambda)h_n |^{1/n} > 1 . $$
\endproclaim
\demo{Proof}
Both $Q_n$ and $r_n$ satisfy the same recurrence relation \thetag{6}.
Hence
$$   Q_{n-1}r_n - Q_nr_{n-1} = -a_1a_2\cdots a_{n-1} , $$
or equivalently
$$   Q_{n-1}h_{n-1} \frac{r_n}{h_na_1\cdots a_n} a_n \frac{h_n}{h_{n-1}}
 - Q_nh_n \frac{r_{n-1}}{h_{n-1}a_1\cdots a_{n-1}} \frac{h_{n-1}}{h_n} = -1. $$
From theorem 1 we get
$$  \left| \frac{r_k}{h_ka_1\cdots a_k} \right| \leq C q^k, \qquad q < 1, $$
so the sequence $Q_nh_n$ can not be majorized by a geometric sequence
$p^n$ with $p < 1/q$, i.e., for any $p < 1/q$ and any positive constant $C$
the inequality
$$   | Q_n(\lambda) h_n | \leq C p^n  $$
is not satisfied for an infinite number of indices $n$. For a subsequence
$\Lambda \subset {\Bbb N}$ we thus have
$$    |Q_n(\lambda)h_n | \geq C p^n , \qquad n \in \Lambda $$
and consequently if we choose $p$ such that $1/q > p > 1$ we have
$$  \limsup_{n \to \infty} |Q_n(\lambda)h_n|^{1/n} \geq p > 1 , $$
giving the required result. \qed
\enddemo

\noindent{\bf Remark 2:} In the symmetric case this corollary gives us
the well known
characterization of $\Omega(A)$ in terms of the orthonormal polynomials
$q_n$: if $\lambda \in \Omega(A)$ then $\limsup_{n \to \infty}
|q_n(\lambda)|^{1/n} > 1$ (\cite{11}).
\medskip

The combination of these two corollaries gives
\proclaim{Theorem 2}
if $A$ is bounded and $0 <  |\alpha_k|, |\gamma_k| \leq C_1 < \infty$
for some constant $C_1$, then for any $\lambda \in \Omega(A)$
there exists a subsequence of the Pad\'e approximants $\pi_n =
P_n(\lambda)/Q_n(\lambda)$ which converges with a geometric rate to
$\phi(\lambda)$.
\endproclaim
\demo{Proof}
From corollary 1 we have
$$  \lim_{n \to \infty} Q_nh_n \phi(\lambda) - P_n h_n = 0 . $$
Let $\Lambda\subset {\Bbb N}$ be such that
$$   \lim_{n \to \infty} |Q_n(\lambda) h_n |^{1/n} > 1, $$
then $|Q_n(\lambda) h_n| \geq C p^n$ for some $p > 1$. Hence
$$   \left| \phi(\lambda) - \frac{P_n(\lambda)}{Q_n(\lambda)} \right|
  \leq C \left( \frac{q}{p} \right)^n, $$
which proves the theorem. \qed
\enddemo

\noindent{\bf Remark 3:} In the symmetric case this statement seems to be
unnoticed.  Indeed, suppose $A$ is in the form \thetag{1}
in an orthonormal basis, and consider the so-called asymptotically
periodic case (or the corresponding limit-periodic continued fraction)
where for some $N$ we have
$$   \lim_{k \to \infty} \alpha_{kN+i} = \delta_i, \qquad
    i=0,1,2, \ldots N-1, $$
and similarly for $\beta_{kN+i}$ and $\gamma_{kN+i}$, then it is known
that the poles of the Pad\'e approximants $\pi_n = P_n/Q_n$ are
essentially concentrated on a system of $N$ closed intervals (\cite{4},
\cite{7}). The convergence
of $\pi_n$ however depends on the behaviour of the so-called spurious
poles. Nevertheless our theorem shows that for any $\lambda$ in the
exterior of the spectrum some subsequence of $\pi_n$ converges
with a geometric rate.
\medskip

Finally we note that from theorem 1 we get a connection between the fact that
$\lambda \in \Omega(A)$ and the convergence of Pad\'e approximants
$\pi_n(\lambda) = P_n(\lambda)/Q_n(\lambda)$. The convergence of
$\pi_n(\lambda)$, without further requirements on the rate, does not imply
that $\lambda \in \Omega(A)$. A counterexample is given by the spectral measure
$d\mu(x) = \sqrt{1-x^2}\, dx$ on $[-1,1]$, for which $P_n(x)= 2^{-n+1}
U_{n-1}(x)$ and $Q_n(x) = 2^{-n} U_n(x)$, where $U_n(x)$ are the Chebyshev
polynomials of the second kind. At $\lambda = 1 \in \sigma(A)$ we have
$\pi_n(1) = 2n/(n+1)$ which converges, but without geometric rate.
The convergence of some subsequence with a geometric
rate is also not sufficient to imply that $\lambda \in \Omega(A)$.
A counterexample is given by any symmetric measure on $[-1,1]$, because
then at $\lambda = 0 \in \sigma(A)$ we always have $P_{2n}(\lambda) = 0$.

The following question is thus of interest: does the convergence of the whole
sequence $\pi_n(\lambda)$ with a geometric rate imply that $\lambda \in
\Omega(A)$?

\enddocument